\documentclass[12pt]{article}
\usepackage{amssymb}
 
\begin{document}

\newtheorem{theorem}{Theorem}
\newcommand{\p}{\partial}
\newcommand{\om}{\omega}
\newcommand{\Om}{\Omega}
\renewcommand{\phi}{\varphi}
\newcommand{\e}{\epsilon}
\renewcommand{\a}{\alpha}
\renewcommand{\b}{\beta}
\newcommand{\N}{{\mathbb N}}
\newcommand{\R}{{\mathbb R}}
\newcommand{\EX}{\mathbb{E}}
\newcommand{\PX}{\mathbb{P}}
 
\title{Ergodicity of  Stochastically Forced Large 
Scale Geophysical Flows  }

\author{Jinqiao Duan$^1$ and Beniamin Goldys$^2$\\
\\ 
1. Department of Applied Mathematics \\
Illinois Institute of Technology \\
Chicago, IL 60616, USA.  \\ 
E-mail: duan@iit.edu\\
\\
2. School of Mathematics \\
The University of New South Wales \\
Sydney 2052, Australia}

\date{March 12, 2001 }

\maketitle

\begin{abstract}

We investigate the  ergodicity of   2D large scale
quasigeostrophic flows under
random wind forcing. We show that  
the quasigeostrophic flows are ergodic  under suitable conditions on the 
random forcing and on the fluid domain, and under no restrictions
on viscosity, Ekman constant
or Coriolis parameter. When these
conditions  are satisfied,  then for any observable of
the   quasigeostrophic flows, its time average   
approximates the statistical ensemble average, as long as  the time 
interval is sufficiently long.

\bigskip


\bigskip


\end{abstract}

\newpage

\section{Introduction}

The models for geophysical flows are usually  very complicated. 
Simplified  models have been developed to investigate
 the basic key features of large scale phenomena. These models  
filter out undesired high frequency oscillations in geophysical flows
and are derived at asymptotically high rotation rate or small Rossby
number. 

An important example of such a geophysical flow model is
the quasigeostrophic flow model  \cite{Pedlosky}  

\begin{eqnarray*}
   \Delta \psi_t + J(\psi, \Delta \psi ) + \beta \psi_x
   =\nu \Delta^2 \psi - r \Delta \psi  + \mbox{wind forcing} \; ,   
\label{oldeqn}
\end{eqnarray*}
where $\psi(x,y,t)$ is the stream function,  
$\beta \geq 0$ is the meridional gradient of the
Coriolis parameter, $\nu >0$ is the viscous dissipation constant,
$r>0$ is the Ekman dissipation constant.
Moreover, $J(f, g)=f_x g_y -f_y g_x$ denotes the Jacobian operator.

The   quasigeostrophic equation   has been
derived as an approximation of the rotating shallow water equations by
the conventional asymptotic expansion in small Rossby number
\cite{Pedlosky}.  
Recently,  the randomly forced
quasigeostrophic flow model has been used to study
various phenomena in geophysical flows   under
uncertain wind forcing 
 \cite{Holloway,Muller,Samelson,Griffa,DelSole-Farrell}.

Introducing (relative) vorticity $\om (x,y,t) = \Delta \psi(x, y,t)$, 
the quasigeostrophic
equation   can be written as
\begin{equation}
       \om_t + J(\psi, \om ) + \beta \psi_x
        =\nu \Delta \om - r \om  +   \mbox{wind forcing}  \; ,
 \label{eqn}
\end{equation}
where $(x, y) \in D$ and $D \subset R^2$ denotes a   bounded domain
with sufficiently regular boundary.  Potential vorticity is defined
as $\om+ \beta y$. The
boundary conditions are no normal flow  ($\psi=0$) 
and free-slip ($\om=0$)
on $\partial D$ as   in  Pedlosky (\cite{Ped96}, p.34) or
in Dymnikov and Kazantsev \cite{Dymnikov}:
\begin{eqnarray}
 \psi = \om =  0  \quad \mbox{on} \; \partial D \; . \label{BC} 
\end{eqnarray} 
An appropriate initial condition  $\om(0)$ is also imposed.
We note that the Poincar\'e inequality holds with these boundary
conditions.

An invariant measure for  stochastic systems is like a
``statistical steady state" and is a part of the asymptotic 
permanent regime
of the system \cite{Arn98}. When there is only one invariant measure
for the   quasigeostrophic flows modeled by (\ref{QG}),
we have the so-called ergodic principle, i.e., for any observable of
the   quasigeostrophic flows, its time average on $[0, T]$
approaches the statistical ensemble average, as $T$ goes to infinity.

We will investigate the existence and uniqueness of invariant measures
for quasigeostrophic flows. After review mathematical setup
in \S 2, we study exixtence and uniqueness of invariant measures
in \S 3 and \S 4, respectively.  Finally, we summarize  our results
in \S 5.

\section{Mathematical Setup}

In the following we use the abbreviations $H=L^2(D)$, 
$H^k_0=H^k_0(D)$, $H^k =H^k(D)$, $ 0 < k < \infty$, for the standard Sobolev
spaces.  Let $<\cdot, \cdot>$ and $\| \cdot \| $
denote the standard scalar product and norm in $L^2$, respectively.
Moreover, the norms for $H^k_0$  are denoted by 
$\| \cdot \|_{H^k}$. 
Due to the Poincar\'e inequality  \cite{Evans}, 
$\| \Delta \phi \|$
is an equivalent norm for  $H^2_0$.
 It is well-known that the linear operator 
$$
A = \nu \Delta : H \to H
$$ 
with domain $D(A) = H^2 \cap H_0^1$ is self-adjoint.  Note
that $A$ generates a strongly continuous, and in fact,
an analytic semigroup $S(t)$ on $L^2$
(\cite{Pazy}).  The spectrum of $ A$ consists of eigenvalues $0 >
\lambda_1 > \lambda_2 \ge \lambda_3 \ge \ldots$ with corresponding
normalized eigenfunctions $e_1$, $e_2$, $\ldots$. The set of
these eigenfunctions is complete in $L^2$. For example, for the square
domain $D = (0,1) \times (0,1)$ the eigenvalues are given by
$-\nu(m^2+n^2)\pi^2$ for positive integers $m,n$, and the associated
eigenfunctions are suitable multiples of $\sin(m\pi x) \sin(n\pi y)$.

We define the nonlinear operator $F$ by 
$$
F(\om) = -r \om -\beta\psi_x -J(\psi, \om ), 
$$
then (\ref{eqn}) can be rewritten as
the abstract evolution equation together with initial condition
\begin{eqnarray}
    d\om & = & (A \om + F(\om) )dt +    \sqrt{Q} dW \; ,  
  		\label{QG}     \\
    & & \om (0)  \;\; \mbox{given}  , 
    		\label{newIC}
\end{eqnarray}
where $W(x,y,t)$ is a  Wiener process defined on a probability space 
$(\Omega, {\cal F}, \mathbb{P})$.  The covariance operator
$Q: H \rightarrow H$   for this Wiener process is a 
nonnegative and symmetric linear continuous operator to be specified below.
The term  with Ito derivative,  $\sqrt{Q} dW$, is a model
for the white-in-time noise representing the random wind forcing.
This equation can be rewritten in the mild (integral) form
\begin{eqnarray}
     \om(t) = S(t) \om(0) + \int_0^t S(t-s)F(\om (s)  )ds +   Z(t) \; ,   
  		\label{mild}  
\end{eqnarray}
where  $Z(t)$ is the 
stochastic convolution
\begin{eqnarray}
   Z(t)  =  \int_0^t S(t-s)\sqrt{Q}dW(s) \; , \quad t>0 \; .
  \label{conv}
\end{eqnarray} 
In fact,  $Z(t)$ is an Ornstein-Uhlenbeck process and it is the
solution of the  linearized version of the above equation (\ref{QG}):
\begin{eqnarray}
    d Z & = &  A Z  \;  dt +  \sqrt{Q} dW \; .  
  		\label{linear}    
\end{eqnarray}

In this paper,  we always assume that the covariance operator
$Q$  for the Wiener process $W(t)$
is of trace class, i.e., $\mbox{Trace} \; Q < +\infty$. 
Thus we only consider noise that is white in time but colored in space.  
Then the
stochastic convolution $Z(t)$ has a continuous version
with values in $H=L^2(D)$; see Theorem 5.14 in \cite{DaPrato}.
 
We can specifically define an appropriate class of Wiener processes $W(t)$
satisfying the above condition.  Let
$\beta_k(t)$, for positive integer  $k$, denote a family of independent real-valued
Brownian motions.  Furthermore, choose positive constant  $\a_k$ 
  such that
\begin{eqnarray}
  \sum_{k=1}^\infty \frac{\a_k^2}{|\lambda_k|^{1-\gamma}} < \infty
\label{sum}
\end{eqnarray}
for some $0 < \gamma < 1$.   
Then we  define the white noise by
 
\begin{eqnarray}
  \sqrt{Q} \dot{W}(t) := \sum_{k=1}^\infty  \a_k  \dot{\beta}_k(t) e_k \; ,
  \quad t \ge 0 \; .
  \label{Wiener}
\end{eqnarray}
Note that the eigenvalues $\lambda_k$
for the operator $A$ behave like $k$ in two dimensions and also note that
the Riemann zeta function $\zeta (s) = \sum_{k=1}^\infty \frac1{k^s}$
is well-defined for $s>1$. We see that the condition (\ref{sum})
is satisfied when $k^{-\frac12} \leq \a_k \leq k^{-\frac38}$.

We further assume that  
$$
\kappa (D) = \inf_{0<\rho < diam(D)} \inf_{(x,y) \in D}
		\frac{meas(D \cap B(x,y; \rho))}{\rho^2} > 0,
$$
where $	diam(D) $ is the diameter of $D$ (the least upper bound
of two-point distances in $D$),   $ meas(\cdot)$
denotes the Lebesgue measure, and $B(x,y; \rho)$
is the open disk centered at $(x,y)$ and with radius $\rho$. 
We also assume that
the eigenfunctions $e_k$ satisfy
\[
e_k \in C_0(\bar{D}), \quad |e_k(x,y)| \leq C,
\]
\[
|\p_x e_k(x,y)|,  \quad |\p_y e_k(x,y)| \leq C \sqrt{|\lambda_k|},
\]
for $(x,y) \in D$,  positive integer $k$, and some constant $C>0$. 
For the square
domain $D = (0,1) \times (0,1)$, these conditions are
all satisfied.
Then, according to Theorem 5.2.9 in \cite{DaPrato2}, the
stochastic convolution $Z(t)$   
has a continuous version with values in $L^2(D)$. 
(Actually, in this case, $Z(t)$ is   in $C_0(D)$,
the Banach space of continuous functions satisfying the zero Dirichlet
boundary condition  on $D$).   For  this Wiener process $W(t)$ 
in (\ref{Wiener}),  the stochastic convolution $Z(t)$ is
\begin{eqnarray}
  Z(t) = \sum_{k=1}^\infty  \a_k  e_k 
     \int_0^t e^{-\lambda_k(t-s)} d \beta_k(s)   \; ,
  \quad t \ge 0 \; .    
\end{eqnarray}

As shown in  \cite{Brannan}, 
for every initial condition $\om(0) \in L^2(D)$, there exists a
unique global mild solution $\om(x,y,t)$ of the quasigeostrophic flow model
(\ref{QG}).
   This solution is  in  $C([0, T]; L^2(D))$ for every $T>0$.

\section{Existence of an Invariant Measure}

Now we consider invariant measure for the quasigeostrophic
flow model (\ref{QG}).  
For the rest of the paper, we   denote $\om(t; x)$ as the solution
of the quasigeostrophic flow  model with {\em initial condition} 
(not the spatial point)
$x \in H$.

We introduce the usual notations.
The Markovian transition semigroup is 
$$
  (P_t g )(x) = \mathbb{E} [g(\om(t; x))],
$$
for $ g \in B_b(H)$, the space of bounded Borel measurable functions. 
Hereafter $\mathbb E$ is the expectation.
The transition probability is
$$
P_t(x, \Gamma) = \PX (\om(t; x) \in \Gamma),
$$
for $x \in H$ and  $\Gamma \in {\cal B}(H)$, the $\sigma$-algebra of Borel sets
in $H$.

A probability measure $\mu$ on $(H, {\cal B}(H))$ is called invariant
if 
$$
\int g d \mu = \int P_t g d \mu
$$
for any $t>0$ and $g \in B_b(H)$, or, equivalently,
$$
\int_H P_t(x, \Gamma) d \mu = \mu (\Gamma),
$$
for any $t>0$, $x \in H$ and $\Gamma \in {\cal B}(H)$.

The existence of an invariant measure  for the quasigeostrophic
flow model (\ref{QG}) follows from
a tightness or , equivalently, a compactness argument \cite{Stroock}.
If the mean-square norm of the solution is bounded for all time
$t>0$ and for all initial data, then by the Chebyshev inequality, 
the solution is bounded in probability, which further implies that
the family of measures on $(H, {\cal B}(H))$
$$
\frac1{T} \int_0^T P_t(x, \cdot) dt,\;\;  T\geq 1,
$$
is tight for some $x \in  H$; see \cite{DaPrato2}, page 89-90.
Thus by Corollary 3.1.2 in 
\cite{DaPrato2}, there exists an invariant measure for the quasigeostrophic
flow model (\ref{QG}). So in the rest of this section, we estimate
the mean-square norm  $ \EX \|\om(t)\|^2$.

We assume that
\begin{eqnarray}
  \int_0^{\infty} \|S(r)\sqrt{Q}\|^2_{HS} dr < +\infty,  
\label{trace} 
\end{eqnarray}
where $\| \cdot \|_{HS}$ is the Hilbert-Schmidt norm.
We rewrite (\ref{mild}) as
\begin{eqnarray}
     \om(t) = Y(t) +   Z(t) \; ,  
\end{eqnarray}
where
$$
Y(t) = S(t) x + \int_0^t S(t-s)F(\om (s)  )ds ,
$$
with initial data $\om(0) =x$,
and $Z(t)$ is the Ornstein-Uhlenbeck process  in  (\ref{conv}).

By Corollary 4.14 in \cite{DaPrato}, for any $x \in H$,
\begin{eqnarray}
\sup_{t \geq 0} \EX \|Z(t)\|^2 = 
\sup_{t \geq 0} \EX \int_0^t \|S(r)\sqrt{Q}\|^2_{HS} dr < +\infty.
\label{Z}
\end{eqnarray}

By \cite{Brannan} or follow a Yosida approximation combined with
$L^2$-norm estimate as in Proposition 6.1.6 in \cite{DaPrato2},
we have, for any $x \in H$,
\begin{eqnarray}
\sup_{t \geq 0}  \EX \|Y(t)\|^2   < +\infty. 
\label{Y}
\end{eqnarray}
Note that
\begin{eqnarray*}
\|\om(t)\|^2 = <Y+Z, Y+Z> & =& \|Y\|^2 + 2<Y, Z> + \|Z\|^2  \\
& \leq & \|Y\|^2 + 2\|Y\| \; \| Z \| + \|Z\|^2  \\
& \leq &  2(\|Y\|^2 +\|Z\|^2).
\end{eqnarray*}
Thus, by (\ref{Z}) and (\ref{Y}),
\begin{eqnarray}
\sup_{t \geq 0}  \EX \|\om(t)\|^2   < +\infty.
\end{eqnarray} 
By the argument in the beginning of this section, there exists
at least one invariant measure for the
the quasigeostrophic
flow model (\ref{QG}). We have the main result in this section.

\begin{theorem}  \label{existence}
Assume that $\int_0^{+\infty} \|S(r)\sqrt{Q}\|^2_{HS} dr < +\infty$.
Then there exists at least one invariant probability measure for  
the quasigeostrophic
flow model (\ref{QG}) in the space $L^2(D)$ of square-integrable 
  vorticities.
\end{theorem}

\section{Uniqueness of an Invariant Measure}

Now we consider the uniqueness of invariant measure 
for  the quasigeostrophic 
flow model (\ref{QG}). As we know in Chapter 4 in
\cite{DaPrato2}, the uniqueness of invariant measure
is a consequence of regularity of the transition semigroup $P_t$,
by the Doob's Theorem. Due to Khasminskii's Theorem,
strong Feller and irreducibility properties imply the regularity.
So we now try to prove the strong Feller and irreducibility properties
for  the transition semigroup $P_t$.

Strong Feller property means that
for every $g(x)$ in $B_b(H)$, the space of bounded Borel 
{\em measurable} functions
on $H$, $P_t g(x)$ is in $C_b(H)$, the space of bounded   
{\em continuous} functions on $H$.

Irreducibility property means that
for every Borel set in $H$, i.e., for every $\Gamma$ in ${\cal B}(H)$,
$P_t(x,\Gamma)$ is positive for any $x \in H$ and $t>0$.

\bigskip
\bigskip

\noindent {\bf  Strong Feller Property}

We first consider strong Feller property.
Note that (\cite{DaPrato}, p.119)
$$
\mbox{Trace} \; \int_0^t S(r)Q S^*(r)dr = \int_0^{t} \|S(r)\sqrt{Q}\|^2_{HS} dr.
$$
So  the condition for the existence of invariant 
measures in Theorem \ref{existence},
i.e.,
 $ \int_0^{+\infty} \|S(r)\sqrt{Q}\|^2_{HS} dr < +\infty$,
implies that the  linear integral operator $Q_t: H \rightarrow H$,
\begin{eqnarray}
 Q_t  x\; :=  \int_0^t S(r)Q S^*(r)x dr, \;\; x \in H,
\label{qt}
\end{eqnarray}
 is of trace class for any $t>0$.

We further assume that 
\begin{eqnarray}
 \mbox{Image} \; S(t)   \subset  \mbox{Image} \; Q_t ^{\frac12}.
\label{image}
\end{eqnarray}
Then follow a similar argument as in the proofs of  Theorem 7.2.4 in
\cite{DaPrato2} and of  Theorem 3.1  in \cite{Ferrario}, we conclude that
$P_t$, $t>0$,  is a strong Feller semigroup.

\bigskip
\bigskip

\noindent  {\bf  Irreducibility  Property }

Now we   consider irreducibility property.
We further assume that the covariance operator $Q$ is
one-to-one (or injective), i.e., the kernel $\ker   Q = \{ 0 \}$.
Then, as in the proof of Theorem 7.4.2 in
\cite{DaPrato2} and of  Theorem 3.1  in \cite{Ferrario},  
$P_t$, $t>0$,  is irreducible.

\bigskip
\bigskip

Thus, with the strong Feller and irreducibility properties proved above,
using Doob's Theorem 4.2.1 in \cite{DaPrato2},
there exists a unique invariant measure $\mu$
on $(H, {\cal B} (H))$, and all other transition probability measures
$P_t(x, \cdot)$, $x \in H$,  approach this unique invariant measure $\mu$
as time goes to infinity.

Therefore, we have the following main theorem in this section.

\begin{theorem} \label{main}
Assume that  

(i) $ \int_0^{+\infty} \|S(r)\sqrt{Q}\|^2_{HS} dr < +\infty$,

(ii) $\mbox{Image} \; S(t)   \subset  \mbox{Image} \; Q_t ^{\frac12}$, 
where $Q_t$ is defined in (\ref{qt}), and

(iii) The covariance operator $Q:\; L^2(D) \rightarrow L^2(D)$ 
is one-to-one.

Then  

(A) There exists a unique invariant probability measure $\mu$    for
the quasigeostrophic
flow system (\ref{QG}) in the  
space  $L^2(D)$  of square-integrable   vorticities;

(B) Moreover,  for any $\om \in  L^2(D)$, 
the transition probability measures
 $P_t(\om, \cdot)$  approach the unique invariant  
probability measure $\mu$.  Namely,   for any $\Gamma \in {\cal B}(H)$,
$$
\lim_{t\rightarrow +\infty} P_t(\om, \Gamma) = \mu (\Gamma);  
$$ and

(C) Quasigeostrophic flow system (\ref{QG})  is ergodic,  namely,
$$
\lim_{T\rightarrow +\infty} \int_0^T g(\om(t)) dt
  = \int_{L^2} g d \mu,  \;\;  \PX - a.s.
$$
for all solution $\om(t)$ with initial date in $ L^2(D)$ and
all Borel measurable function $g: L^2(D) \rightarrow  R$
such that $\int_{L^2(D)} \|g\| d \mu < \infty$.
\end{theorem}

The ergodicity in Part (C) above is a consequence of the uniqueness of 
the invariant measure $\mu$; see Theorem 3.2.6 in
\cite{DaPrato2}.

\section{Summary}

In this paper, we have studied ergodicity of
large scale quasigeostrophic flows
under random wind forcing.  We have shown that  
the quasigeostrophic flows are ergodic   under suitable conditions on the 
random forcing and on the fluid domain, and under no restrictions
on viscosity, Ekman constant
or Coriolis parameter.  When these
conditions  are satisfied,  then for any observable of
the   quasigeostrophic flows, its time average   
approximates the statistical ensemble average, as long as the time 
interval is sufficiently long.

There is recent work on random dynamical attractors for the 
 quasigeostrophic flow model
by Duan et al. \cite{Duan-Kloeden-Schmalfuss}.
A consequence of that work implies that, when viscosity is sufficiently
 large
and  when the trace of the covariance operator for the Wiener process
is sufficiently small,  then all quasigeostrophic motions 
approach a  point random attractor exponentially fast as
time goes to infinity. This is a very rare case.
This  point random attractor   corresponds to a
unique invariant Dirac measure, i.e., the supporting point of the
Dirac measure is the global (point) attractor, and thus
under these conditions, quasigeostrophic flows are also ergodic.
 These conditions are different from the ergodic
conditions in the current paper.  For example, in the current paper,
we do not
impose any condition  on viscosity,  or on the size of
the trace of   the covariance operator for the Wiener process.

\bigskip 
\bigskip 
\bigskip 

{\bf Acknowledgement.} A part of this work was done while J. Duan was
visiting the University of New South Wales, Australia.  
This work was partly supported by the Australia Research Council   and
by the NSF Grant 
DMS-9973204.


\begin{thebibliography}{60}
 
 
\bibitem{Arn98}
L.~Arnold.
\newblock {\em {R}andom {D}ynamical {S}ystems}.
\newblock Springer--Verlag, Berlin, 1998.

 
\bibitem{Brannan}
J.~R. Brannan, J. Duan, and T.~Wanner.
\newblock Dissipative quasigeostrophic dynamics under random forcing.
\newblock {\em J. Math. Anal. Appl.}, 228:221--233, 1998.
 


\bibitem{DaPrato} G. Da Prato and J. Zabczyk, {\em Stochastic Equations 
in Infinite Dimensions}, Cambridge University Press, 1992.


\bibitem{DaPrato2} G. Da Prato and J. Zabczyk, {\em Ergodicity for  
Infinite Dimensional Systems}, Cambridge University Press, 1996.

 


\bibitem{DelSole-Farrell} T. DelSole and B. F. Farrell,
A stochastically excited linear system as a model for quasigeostrophic 
turbulence: Analytic results for one- and two-layer fluids,
{\em J. Atmos. Sci.} {\bf 52} (1995) 2531-2547.
 
 
 
\bibitem{Duan-Kloeden-Schmalfuss} 
\newblock J. Duan, P. E. Kloeden and B. Schmalfu{\ss}. 
\newblock  Exponential Stability of the Quasigeostropic Flows under Random Perturbation,  {\em Progress in Probability}, in press, 2000.


\bibitem{DuanSchm} 
J.~Duan and B.~Schmalfu{\ss}.
\newblock   The 3D Quasigeostrophic Equation under Random
Perturbation.
\newblock {\em     },   submitted, 2000.

\bibitem{Dymnikov} V.  Dymnikov and E. Kazantsev,
On the genetic ``memory" of chaotic attractor of the
barotropic ocean model,
In Proceedings of the third bilateral conference 
``Predictability of atmospheric and oceanic circulations " 
of the French-Russian A.M.Liapunov Institute in Computer Science
 and Applied Mathematics (INRIA - Moscow State University). Nancy, 
April, 1997.Edition MSU, 1997, pp. 25-36.
  
\bibitem{Evans} L. C. Evans,
{\em Partial Differential Equations}, 
Amer. Math. Soc., Providence, 1998.

\bibitem{Ferrario} B. Ferrario, 
Ergodic results for stochastic Navier-Stokes equation,
{\em Stoch. Stochastic reports} {\bf 60} (1997), 271-288.

\bibitem{Ferrario2} B. Ferrario,
Stochastic Navier-Stokes equations: analysis of the noise to have a unique
invariant measure. {\em Ann. Mat. Pura Appl.} (4) 177 (1999), 331--347.


 

\bibitem{F}
F. Flandoli,
Dissipativity and invariant measures for stochastic 
 Navier-Stokes equations,
{\em NoDEA.} {\bf 1 } (1994),  403-423. 


\bibitem{FM}
F. Flandoli and B. Maslowski,
Ergodicity of the 2D Navier-Stokes equation under
random perturbations,
{\em Commun. Math. Phys.} {\bf 171} (1995), 119-141.



\bibitem{Griffa} A. Griffa and S. Castellari,
Nonlinear general circulation of an ocean model driven by
wind with a stochastic component,
{\em J. Marine Research} {\bf 49} (1991), 53-73.

 


\bibitem{Holloway} G. Holloway, Ocean circulation: Flow in probability under
statistical dynamical forcing,
in {\em Stochastic Models in Geosystems},
S. Molchanov and W. Woyczynski (eds.), Springer, 1996.

  

\bibitem{Muller} P. M\"uller, Stochastic forcing of quasi-geostrophic eddies,
in {\em Stochastic Modelling in Physical Oceanography},
R. J. Adler, P. M\"uller and B. Rozovskii (eds.), Birkh\"auser, 1996.


\bibitem{Pazy} A. Pazy, {\it Semigroups of linear operators and
applications to partial differential equations}, Springer-Verlag,  1983.

 
\bibitem{Pedlosky} J. Pedlosky, {\it Geophysical Fluid Dynamics},
Springer-Verlag, 2nd edition, 1987.


\bibitem{Ped96}  J.~Pedlosky,
 {\em Ocean Circulation Theory}.
Springer--Verlag,  Berlin, 1996.




\bibitem{Samelson} R. M. Samelson, Stochastically forced current 
fluctuations in vertical shear and over topography,
{\em J. Geophys. Res.} {\bf 94} (1989) 8207-8215.

 
 
 

 \bibitem{Stroock} D. W. Stroock,
{\em Probability Theory: An Analytic View},
Cambridge University Press, Cambridge, 1993.


\bibitem{Vishik} M. I. Vishik and A. V. Fursikov,
{\em Mathematical Problems of Statistical  Hydrodynamics},
Kluwer Academic Publishers, Boston, 1988. 
 
\end{thebibliography}
\end{document}